\documentclass{amsart}

\usepackage{amsopn}
\usepackage{amssymb}
\usepackage{amscd}

\def\Vs{{\scriptscriptstyle V}}

\newtheorem{theorem}{Theorem}[section]
\newtheorem{lemma}[theorem]{Lemma}
\newtheorem{proposition}[theorem]{Proposition}

\newtheorem*{1.1.1}{\textbf{Theorem 1.1.1}}
\newtheorem*{1.1.2}{\textbf{Proposition 1.1.2}}
\newtheorem*{1.1.3}{\textbf{Proposition 1.1.3}}
\newtheorem*{1.1.4}{\textbf{Corollary 1.1.4}}
\newtheorem*{1.2.1}{\textbf{Theorem 1.2.1}}
\newtheorem*{1.3.1}{\textbf{Theorem 1.3.1}}
\newtheorem*{1.4.1}{\textbf{Proposition 1.4.1}}
\newtheorem*{1.4.2}{\textbf{Theorem 1.4.2}}

\theoremstyle{definition}

\theoremstyle{remark}

\numberwithin{equation}{section}

\begin{document}

\title[Embeddings of proper actions in compactifications]{On embeddings of proper and equicontinuous actions  in zero-dimensional compactifications}

\author[Antonios Manoussos]{Antonios Manoussos}
\address{Fakult\"{a}t f\"{u}r Mathematik, Universit\"{a}t Bielefeld, SFB 701, Postfach 100131, D-33501 Bielefeld, Germany}
\email{amanouss@math.uni-bielefeld.de}
\thanks{This work was partially supported by DFG Forschergruppe ``Spektrale Analysis, asymptotische Verteilungen und
stochastische Dynamik" and SFB 701 "Spektrale Strukturen und Topologische Methoden in der Mathematik" at the University of Bielefeld, Germany.}

\author[Polychronis Strantzalos]{Polychronis Strantzalos}
\address{Department of Mathematics, University of Athens, Panepistimioupolis, GR-157 84, Athens, Greece}
\email{pstrantz@math.uoa.gr}
\thanks{}


\date{}


\begin{abstract}
We provide a tool for studying properly discontinuous actions of non-compact groups on locally compact, connected and paracompact spaces, by embedding such an action
in a suitable zero-dimensional compactification of the underlying space with pleasant properties. Precisely, given such an action $(G,X)$ we construct a
zero-dimensional compactification $\mu X$ of $X$ with the properties: (a) there exists an extension of the action on $\mu X$, (b) if $\mu L\subseteq \mu X\setminus X$
is the set of the limit points of the orbits of the initial action in $\mu X$, then the restricted action $(G,\mu X\setminus \mu L)$ remains properly discontinuous,
is indivisible and equicontinuous with respect to the uniformity induced on $\mu X\setminus \mu L$ by that of $\mu X$, and (c) $\mu X$ is the maximal among the
zero-dimensional compactifications of $X$ with these properties. Proper actions are usually embedded in the end point compactification  $\varepsilon X$ of $X$, in
order to obtain topological invariants concerning the cardinality of the space of the ends of $X$, provided that $X$ has an additional ``nice" property of rather
local character (``property Z", i.e., every compact subset of $X$ is contained in a compact and connected one). If the considered space has this property, our new
compactification coincides with the end point one. On the other hand, we give an example of a space not having the ``property Z" for which our compactification is
different from the end point compactification. As an application, we show that the invariant concerning the cardinality of the ends of $X$ holds also for a class of
actions strictly containing the properly discontinuous ones and for spaces not necessarily having ``property Z".
\end{abstract}

\maketitle

\section*{Introduction}
The end point compactification of a locally compact space has been proved fruitful for the study of the space in the topological framework, including proper actions.
One reason for this is that we have a ``clear view" of the embedded space in such a compactification, contrary to the situation when, for example, the Stone-\v{C}ech
compactification is consider instead. Actually, the end point compactification  is the quotient space of the Stone-\v{C}ech compactification with respect to the
equivalence relation whose equivalence classes are the singletons of $X$ and the connected components of $\beta X\setminus X$.

Our purpose in this paper is to provide an equivariant and analogously useful notion corresponding to the end point compactification in order to have a ``clear view"
of the embedded proper action. By saying a ``clear view" of the embedded action we mean that the embedded action has at least the three properties that follow:

Let $(G,X)$ be the initial proper action, $(G,Y)$ the extended action in a zero-dimensional compactification $Y$ of $X$ and $L$ be the set of the limit points of the
orbits of the initial action in $Y$ (i.e., the cluster points of the nets $\{ g_{i}x\}$, for all nets $\{ g_{i}\}$ divergent in $G$, and $x\in X$), then the maximal
invariant subspace where the extended action can be proper is, obviously, $Y\setminus L\supseteq X$. So, the required properties are: The action $(G,Y\setminus L)$

(a) remains proper

(b) is equicontinuous with respect to the uniformity induced on $Y\setminus L$ by that of $Y$, and

(c) is indivisible, (i.e., if $\lim g_{i}y_{0}=e\in L$ for some $y_{0}\in Y\setminus L$ then $\lim g_{i}y=e$ for every $y\in Y\setminus L$).

In this direction the main results of the paper at hand are:

(1) \hspace{0.01cm} If $X$ is a locally compact, connected and paracompact space and $G$ is a non-compact group acting properly discontinuous on $X$, there always
exists a zero-dimensional compactification $\mu X$ of $X$ which is the maximal (in the ordering of zero-dimensional compactifications of $X$) that satisfies the
properties: (a) the initial action can be extended on $\mu X$, and (b) if $\mu L$ denotes the set of the limit points of the orbits of the initial action in $\mu X$,
the restricted action $(G,\mu X\setminus \mu L)$ remains proper, is equicontinuous with respect to the uniformity induced on $\mu X\setminus \mu L$ by that of $\mu X$
and indivisible as embedded in the action $(G,\mu X)$ (Theorem 6.2).

(2) \hspace{0.01cm} $\mu L$ consists of at most two or infinitely many points (Theorem 6.3).

(3) \hspace{0.01cm} If $X$ has the \textit{``property Z"}, i.e., every compact subset of $X$ is contained in a compact and connected one (for example if $X$ is
locally compact, connected and locally connected) then $\mu X$ coincides with the end point compactification $\varepsilon X$ of $X$ (Corollary in Section 6).

The proof of the results stated above relies on a \textit{new construction}: The action $(G,\mu X)$ is obtained by taking the initial action as an equivariant inverse
limit of properly discontinuous $G$-actions on polyhedra, which are constructed via $G$-invariant locally finite open coverings of $X$, generated by locally finite
coverings of (always existing) suitable fundamental sets of the initial action (cf. Section 3).

As an application of these results we prove in Theorem 7.1 that the invariant concerning the cardinality of the ends for proper actions of non-compact groups on
locally compact and connected spaces with the ``property Z", holds also for proper actions on spaces not necessarily satisfying this property:

If either $G_{0}$, the connected component of the neutral element of $G$, is non-compact, or $G_{0}$ is compact and $G/G_{0}$ contains an infinite discrete subgroup,
then $X$ has at most two or infinitely many ends.

Moreover, in Section 2 we give an example of a properly discontinuous action $(G,X)$, where $G$ is a non-compact group and $X$ is a locally compact, connected and
paracompact space not satisfying the ``property Z" such that $\mu X$ does not coincide with the end point compactification of $X$: we show that the sets of the limit
points of the actions $(G,X)$ and $(G,\varepsilon X\setminus \varepsilon L)$ in $\varepsilon X$ coincide, but the action $(G,\varepsilon  X\setminus \varepsilon L)$
is neither proper nor equicontinuous with respect to the uniformity induced on $\varepsilon X\setminus \varepsilon  L$ by that of $\varepsilon X$.

Properties (a) and (b) in (1) above have been used already especially concerning embeddings in the end point compactification, in order to prove that the existence of
a proper action $(G,X)$ of a non-compact group on a locally compact, connected and paracompact space with the ``property Z" has implications in the structure and the
cardinality of the space of ends of $X$. The following indications trace the known results in this direction:

The first theorem that relates, although indirectly, equicontinuous actions with structural features of spaces, is formulated by Ker\'{e}kj\'{a}rt\'{o} (1934), who
proved that, if the abelian group generated by a homeomorphism of the 2-sphere, $S^{2}$, acts equicontinuously on $S^{2}$ with respect to the metric uniformity of
$S^{2}$ except for a finite number of points, then the number of these exceptional points is at most two. These points can be viewed as the set of the end points of
the maximal subspace of $S^{2}$ on which the above group acts equicontinuously. This result is considerably generalized by Lam in \cite{Lam1} for equicontinuous
actions of non-compact groups on locally compact, connected metric spaces $X$ with respect to uniformities induced, say, by the uniformities of suitable
zero-dimensional compactifications of $X$, i.e., compactifications with zero-dimensional remainder. Roughly speaking it is shown that, if an action $(G,X)$ can be
embedded in an action $(G,Y)$, where $Y$ is a zero-dimensional compactification of $X$, such that

(a) \hspace{0.01cm} there exists a subset $R\supseteq X$ of $Y$ such that $Y\setminus R$ is exactly the non-empty set of the points where the action $(G,Y)$ is not
equicontinuous, and

(b) \hspace{0.01cm} the restricted action $(G,R)$ is indivisible (i.e., if $\lim g_{i}y_{0}=e\in Y\setminus R$ for some $y_{0}\in R$ then $\lim g_{i}y=e$ for every
$y\in R$),

\noindent then $Y\setminus R$ consists of at most two or infinitely many points.

On the other hand, similar results are proved by Abels in \cite{Ab2} for proper actions $(G,X)$, where $G$ was a non-compact topological group and $X$ was a locally
compact and connected space with the ``property Z" (i.e., every compact subset of $X$ is contained in a compact and connected one). The corresponding property in
Lam's work was $X$ to be semicontinuum , which ensured the indivisibility of the equicontinuous action on $R$. In \cite{Ab2} it is considered the end point
compac\-tification, $\varepsilon X$, of $X$, the maximal compactification of it with zero-dimensional remainder, instead of an appropriate zero-dimensional
compactification $Y$ of $X$, and it is proved that such a proper action $(G,X)$ has an extension on $\varepsilon X$. To be more precise, let $\varepsilon L$ denotes
the set of the limit points of the action $(G,X)$ in $\varepsilon X$. Then, it was shown that (a) the action $(G,\varepsilon X\setminus \varepsilon L)$ remains proper
and (b) it is indivisible. Using this embedding, it was shown that $X$ has at most two or infinitely many ends, a remarkable invariant of the proper action $(G,X)$ of
the non-compact group $G$.

The interconnection of the main results in \cite{Lam1} and \cite{Ab2} is explained in \cite{Str1}, where it was shown that, for spaces with the ``property Z", a group
acting equicontinuously in Lam's view may considered as a dense (not necessarily strict) subgroup of a group acting properly as in Abels' view.

\section{Preliminaries}
\noindent\textbf{1.1.} The \textit{Freudenthal} or \textit{end point} compactification $\varepsilon X$ of a locally compact space $X$  may defined as the quotient
space of the Stone-\v{C}ech compactification $\beta X$ of $X$ with respect to the equivalence relation whose equivalence classes are the singletons of $X$ and the
connected components of $\beta X\setminus X$. Recall that the zero-dimensional compactifications of $X$ are ordered with respect to the following ordering: Let $Y$
and $Z$ be two zero-dimensional compactifications of $X$; then $Y\leq Z$ if there exists a surjection from $Z$ onto $Y$ extending the identity map of $X$. Therefore,
the end point compactification is the maximal zero-dimensional compactification of $X$, i.e., for every zero-dimensional compactification $Y$ of $X$ there is a
surjection $p:\varepsilon X\rightarrow Y$ extending the identity map of $X$.

The points of $\varepsilon X\setminus X$ are the \textit{ends} of $X$.

\vspace{0.3cm}

The following theorem, \cite{Ray}, provides an equivalent definition.

\begin{1.1.1} If $Y$ is a compactification of $X$, it is the
end point compactification of $X$ iff $Y\setminus X$ is totally disconnected and does not disconnect $Y$ locally, i.e., given an open (in $Y$) neighborhood $V$ of
$y\in Y\setminus X$ then there is no decomposition of $V\cap X$ into two open disjoint subsets $U_{1},$ $U_{2}$ such that $y\in \overline{U_{1}}\cap
\overline{U_{2}}.$
\end{1.1.1}

The end point compactification has the following useful properties.

\begin{1.1.2} Let $X$ and $Y$ be two locally compact
topological spaces. Then every proper map  $f:X\rightarrow Y$ may extended to a unique map $\varepsilon f: \varepsilon X\rightarrow \varepsilon Y$ that maps ends of
$X$ to ends of $Y$.
\end{1.1.2}

\begin{proof} By the characteristic property of the Stone-\v{C}ech compactification, the map $f:X\rightarrow Y$ has a unique extension $\varepsilon f:\varepsilon
X\rightarrow \varepsilon Y$. The inclusion $\varepsilon f(\varepsilon X\setminus X)\subseteq \varepsilon Y\setminus Y$ follows from the assumption that $f$ is a
proper map.
\end{proof}

\begin{1.1.3} Let $X$ be  a locally compact and connected
space and $Y$ be a zero-dimensional compactification of $X$. Then, whenever a continuous action $(G,X)$ has an extension $(G,Y)$ this extension is continuous.
\end{1.1.3}

\begin{proof} It suffices to show the continuity of the extended action map at the point $(e,z)$, where $e$ is the neutral element of $G$ and $z\in Y$. Let $V$ and
$U$ be two open neighborhoods of $z$ in $Y$ with boundaries in $X$ such that $\overline{V}\subseteq U$. Since the boundaries $\partial U$ and $\partial V$ are compact
subsets of $X$, the set $A=\{g\in G\,|\, g\partial V\subseteq U\, \mbox{and}\, g^{-1}\partial U\subseteq Y \setminus\overline{V}\}$ is an open neighborhood of $e$ in
$G$. We shall show that $g\overline{V}\subseteq U$ for every $g\in A$: The boundary of the set $g\overline{V}\cap (Y\setminus U)$ is contained in $(g\partial V\cap
(Y\setminus U))\cap (g\overline{V}\cap\partial U)$, which is empty by the definition of $A$. Since $Y$ is connected, this implies that $g\overline{V}\cap (Y\setminus
U)$ is either the empty set or coincides with $Y$. The latter is impossible since, choosing a point $x\in\partial V$,  the definition of $A$ implies that $gx\notin
Y\setminus U$. Therefore $g\overline{V}\cap (Y\setminus U)=\emptyset$.
\end{proof}

As an immediate consequence of the above two propositions we state the following

\begin{1.1.4} An action $(G,X)$ of a group $G$ on a locally
compact and connected space $X$ may extended to a unique action on the end point compactification  $\varepsilon X$ of $X$.
\end{1.1.4}

\noindent\textbf{1.2.} The notion of a proper action is given in \cite[III, 4]{Bour1}. Equivalently, an action $(G,X)$ is \textit{proper} if $J(x)$ is the empty set
for every $x\in X$, where
$$
\left.\begin{array}{c} J(x)=\{ y\in X\,|\, \mbox{there exist nets}\,\, \{x_{i}\}\,\, \mbox{in}\,\, X\,\, \mbox{and}\,\, \{g_{i}\}\,\, \mbox{in}\,\, G\,\,
\mbox{with}\,\, g_{i}\rightarrow\infty,\\
\lim x_{i}=x\,\, \mbox{and}\,\, \lim g_{i}x_{i}=y\}.
\end{array}\right.
$$
Here $g_{i}\rightarrow\infty$ means that the net $\{g_{i}\}$ does not have any limit point in $G$.

\vspace{0.3cm}

In the special case where $G$ is locally compact, an action $(G,X)$ is proper iff for every $x,y\in X$ there exist neighborhoods $U_{x}$ and $U_{y}$ of $x$ and $y$,
respectively, such that the set
$$
G(U_{x},U_{y})=\{g\in G\,|\,(gU_{x})\cap U_{y}\neq\emptyset\}
$$
is relatively compact in $G$.

The action is called \textit{properly discontinuous} when
$G(U_{x},U_{y})$ is finite.\\
\\
\textbf{Remark} Let $(G,X)$ be a proper action of a non-compact group $G$ and $(G,\varepsilon X)$ its extension on the end point compactification of $X$. Then, the
set $J(x)$ with respect to the extended action is a non-empty subset of $\varepsilon X\setminus X$ for every $x\in X$. The study of these sets provides useful
information. As an example, we note from \cite{Ab2} the following

\begin{1.2.1} Let $(G,X)$ be a proper action of a non-compact
group $G$ on a locally compact and connected space $X$ with the ``property Z". Then, $X$ has at most two or infinitely many ends. Especially, if $G$ is connected,
then $X$ has at most two ends.
\end{1.2.1}

\noindent\textbf{1.3.} A characteristic and very useful feature of a proper
action is the fundamental set (cf. \cite{Koszul} and \cite{Ab1}).\\

\noindent\textbf{Definition} Given an action $(G,X)$, a subset $F$ of $X$ is a \textit{fundamental set} for the action if $GF=X$ and for every compact subset $K$ of
$X$ the set $\{g\in G\,|\, (gK)\cap F\neq\emptyset\}$ is relatively compact in $G$.

\vspace{0.3cm}

The existence of a fundamental set implies that the action $(G,X)$ is proper, but the converse does not hold, in general. The notion of the fundamental set is
relative to the well known notion of a section but is different in general, in the sense that there are cases where a section is a fundamental set, a fundamental set
fails to be a section and cases where a section fails to be a fundamental set.

\begin{1.3.1}  Let $(G,X)$ be a proper action, where $X$ is a
locally compact, connected and paracompact space. Then, there exist open fundamental sets $F$ and $S$ for $G$ in $X$ such that $\overline{F}\subseteq S$.
\end{1.3.1}

This follows immediately by \cite[Lemma 2, p. 8]{Koszul}, because $X$ is $\sigma$-compact, hence the orbit space  of the action is paracompact.

\vspace{0.5cm}

\noindent\textbf{1.4.} Establishing the notation, we recall\\

\noindent\textbf{Definition} An inverse system $(X_{\lambda },p_{\kappa \lambda },\Lambda )$ consists of a directed set $\Lambda $, a family of topological spaces
$\{X_{\lambda },$ $\lambda \in \Lambda \}$, and continuous mappings $p_{\kappa \lambda }:X_{\lambda }\rightarrow X_{\kappa }$ with the properties that for every
$\kappa ,\lambda ,\mu \in \Lambda$ with $\kappa \leq \lambda $ and $\lambda \leq \mu $ the map $p_{\lambda \lambda }:X_{\lambda }\rightarrow X_{\lambda } $ is the
identity of $X_{\lambda }$, and $p_{\kappa \lambda }\circ p_{\lambda \mu }=$ $p_{\kappa \mu }$. Let $p_{\lambda }:\prod_{\lambda } X_{\lambda }\rightarrow X_{\lambda
}$ be the $\lambda$-projection. The (possibly empty) space
$$
\{x\in \prod_{\lambda } X_{\lambda }\,|\,p_{\kappa }(x)=p_{\kappa \lambda }\circ p_{\lambda }(x)\,\,\mbox{for every}\,\kappa \leq \lambda \}
$$
is called the \textit{inverse limit} of $\{X_{\lambda },$ $\lambda \in \Lambda \}$ and is denoted by $\varprojlim X_{\lambda }.$

\begin{1.4.1}[{\cite[Pr. 2.3, p. 428]{Dug1}}] The sets $\{p_{\lambda
}^{-1}(U)\, |\,\lambda \in \Lambda ,\,U\,\mbox{open}\,\mbox{in}\,X_{\lambda }\}$ form a basis for $\varprojlim X_{\lambda }$.
\end{1.4.1}

The following notion provides an alternative way to describe
locally compact and paracompact spaces using coverings.\\

\noindent\textbf{Definition} Let $X$ be a paracompact space, $(X_{\lambda },p_{\kappa \lambda },\Lambda )$ be an inverse system and
$\{p_{\lambda}\,|\,\lambda\in\Lambda\}$ a family of mappings $p_{\lambda }:X\rightarrow X_{\lambda }$ such that $p_{\kappa }(x)=p_{\kappa \lambda }\circ p_{\lambda
}(x)$ for every $\kappa \leq \lambda$.  We say that the inverse system $(X_{\lambda },p_{\kappa \lambda },\Lambda )$ is a \textit{resolution} of $X$ if the following
conditions hold:

(a) \hspace{0.01cm}  For every covering $\mathcal{U}$ of $X$ that admits a subordinated partition of unity, there exist an index $\lambda \in \Lambda$ and a covering
$\mathcal{U}_{\lambda }$ of $X_{\lambda}$ that also admits a subordinated partition of unity, such that $p_{\lambda}^{-1}(\mathcal{U}_{\lambda })$ refines
$\mathcal{U}$.

(b)  \hspace{0.01cm} For every $\kappa \in \Lambda$ and every covering $\mathcal{U}_{\kappa}$ of $X_{\kappa}$, as above, there exists $\lambda \geq \kappa$ such that
$p_{\kappa \lambda}(X_{\lambda}) \subseteq St(p_{\kappa}(X),\mathcal{U}_{\kappa })$, where
$$
St(B,\mathcal{U})=\bigcup \{U_{i}\,|\, U_{i}\cap B\neq \emptyset,\, U_{i}\in \mathcal{U}\}
$$
is the \textit{star} of $B$ with respect to the covering $\mathcal{U}$.

\begin{1.4.2}[{\cite[Cor. 4, p. 83]{Mar1}}] If the spaces
$X_{\lambda}$ are normal and $X$ is paracompact, then a resolution of $X$ gives $X$ as an inverse limit.
\end{1.4.2}

\section{A counterexample}
Following the notation in the introduction, we now give an example showing that, if the space $X$ does not have the ``property Z", then the action $(G,\varepsilon
X\setminus \varepsilon L)$ is not necessarily proper or equicontinuous.\\
\\
\textbf{2.1.} \textit{The half-open Alexandroff square} $Y$ is the space $[0,1]\times [0,1]\setminus \{ (x,y)\,|\,x=0\, \mbox{and}\, y\in (0,1]\, \mbox{or}\, x=1\,
\mbox{and}\, y\in [0,1)\}$ endowed with the topology $\tau$ defined as follows: A neighborhood basis  of a point $(x,x)\in\Delta=\{ (x,x)\,|\,x\in [0,1]\}$ is
obtained by the intersection of $Y$ with open (in $Y\subseteq\mathbb{R}^{2}$) horizontal stripes less a finite number of vertical lines; a neighborhood basis for the
points $p=(s,t)$ off $\Delta$ is obtained by the intersection of $Y\setminus\Delta$ with open vertical segments centered at $p$ (cf. \cite[Ex. 101, p. 120]{Steen}).
This space is a compact, connected and not locally connected Hausdorff space. Observe that, if $\{(x_{i},y_{i})\}$ is a net converging with respect to the Euclidean
topology on $Y$ to a point $(x,y)$, then this net converges to $(y,y)$ with respect to $\tau$, unless there is an index $i_{0}$ such that $x_{i}=x$ for all $i\geq
i_{0}$, in
which case it converges to $(x,y)$.\\
\\
\noindent\textbf{2.2.} Let $X$ be the subspace $(0,1)\times (0,1)$ of $Y$. This space is locally compact, connected and paracompact, because the closed horizontal
stripes are compact subsets of $X$. It does not have the ``property Z", because every closed horizontal stripe is not contained in a compact and connected subset of
$X$.

The space $Y$ is the end point compactification, $\varepsilon X$, of $X$, and the ends are the points $(x,0)$ for $x\in [0,1)$ and $(x,1)$ for $x\in (0,1]$. In order
to prove this, by Theorem 1.1.1, it is sufficient to verify that the set $Y\setminus X$ is totally disconnected and that every point of it does not disconnect $Y$
locally: For the points of the form $(x,0)$ and $(x,1)$ for $x\in (0,1)$ this follows from the fact that a neighborhood basis of every one of these points consists of
half-open vertical segments which they do not disconnect $Y$ locally. To verify the same for the points $(0,0)$ and $(1,1)$ observe that, if there is a neighborhood
$V$ (in $Y$) for, e.g., $(0,0)$ such that $V\cap X$ is the union of two open sets (in $X$) having $(0,0)$ as a common point of their closures in $Y$, then
they have common interior points.\\
\\
\textbf{2.3.} Next we define a properly discontinuous action of the additive group of the integers $\mathbb{Z}$ on $X$. For convenience, we consider $X$ as
$\mathbb{R}^{2}$ endowed with the topology $\tau$, and we define the action by letting
\begin{center}
$z(x,y)=(x+z,y+z)$ for $z\in\mathbb{Z}$ and $(x,y)\in\mathbb{R}^{2}$.
\end{center}
By Corollary 1.1.4, this action has an extension on $Y=\varepsilon X$. The set $\varepsilon L$, of the limit points of this action, consists of the points $(0,0)$ and
$(1,1)$. The restricted action $(\mathbb{Z},\varepsilon X\setminus \varepsilon L)$ is neither proper nor equicontinuous with respect to the uniformity induced on
$\varepsilon X\setminus \varepsilon L$ by that of $\varepsilon X$. For this, observe that the sequence $\{ (x-n,x)\,|\,n\in \mathbb{N}\}$ converges to the point
$(x,x)$, while the sequence $\{n(x-n,x)=(x,x+n)\}$ converges to an end $e$ that corresponds to the vertical line $\{ (x,y)\,|\,y\in \mathbb{R}\}$. Since $e\in
J((x,x))$, the action $(\mathbb{Z},\varepsilon X\setminus \varepsilon L)$ is not proper. On the other hand, $\lim n(x-n,x)=e$ and $\lim n(x,x)=(1,1)$, therefore this
action is not equicontinuous at $(x,x)$.

\section{The basic construction}
In the sequel we shall proceed to answer the question formulated in the introduction. Our answer will be  based on an inverse system of properly discontinuous actions
on polyhedra, defined from the initial action on $X$. This is achieved using appropriate invariant locally finite coverings of the given space, in order to have the
initial action as an inverse limit of them. To obtain this, it is reasonable to work with invariant coverings of $X$ extending specific coverings of always existing
fundamental sets of the initial action. The construction of this inverse system, which follows, is new and will be given
in several steps:\\
\\
\textbf{3.1.} Let $(G,X)$ be a properly discontinuous action of a non-compact group $G$ on a locally compact, connected and paracompact space $X$. Recall that a
covering $\mathcal{V}$ of $X$ is called a barycentric refinement of a given covering if the covering $\{St(x,\mathcal{V})\,|\,x\in X\}$ refines it, where
$St(x,\mathcal{V})$ has been defined in 1.4. Since $X$ is a locally compact and paracompact space, by \cite[Cor. 7.4, p. 242]{Dug1}, starting with an open covering of
$X$, we can always find an open locally finite barycentric refinement $\mathcal{V}= \{V_{j}\,|\,j\in J\}$ of it consisting
of relatively compact open sets.\\
\\
\textbf{3.2.} Theorem 1.3.1 ensures that there exist open fundamental sets $F$ and $S$ such that $\overline{F} \subset S$. With the previous notation, we can choose
$\mathcal{V}$ such that (a) if $V_{j}$ intersects the boundary of the open fundamental set $S$, then $V_{j}$ does not intersect the open fundamental set $F$, and (b)
if $V_{j}$ does not intersect the boundary of $S$ then either $V_{j}\subseteq S$ or $V_{j}\subseteq X\setminus {\overline S}$. The family $\mathcal{U}=
\{St(x,\mathcal{V})\,|\,x\in F\}$  is an open locally finite covering of $F$ (in $X$); it is also a covering of ${\overline F}$, because if some $V_{j}$ intersects
the boundary of $F$ then intersects $F$, hence is a member of a star of some point of $F$. From (a) and (b) it is easily seen that each member of $\mathcal{U}$
is a subset of $S$.\\
\\
\textbf{3.3.} In the sequel we shall use the following modification of the previous construction, aiming to enrich $\mathcal{U}$ with the property: if $U_{i}\in
\mathcal{U}$ and $gU_{i}\cap F\neq\emptyset$ for some $g\in G$, then $gU_{i}\in \mathcal{U}$. To this end, let $\mathcal{W}=\{ W_{k}\,|\, k\in K\}$ be a locally
finite refinement of $\mathcal{V}$ with the property that the closures of the stars of it are subsets of corresponding stars of $\mathcal{V}$. Now, we observe that
the set $M_{i}=\{g\in G\,|\,gU_{i}\cap {\overline F}\neq\emptyset \}$ is non-empty and finite, because the action is properly discontinuous, $U_{i}$ is relatively
compact and ${\overline F}\subseteq S$ (cf. 1.2 and 1.3).

If $x\in {\overline W_{k}}\subseteq U_{i}$ for some $U_{i}\in \mathcal{U}$, and $gx\in \overline{F}$, then $g\in M_{i}$ which is finite. So, for $x\in X$ we can find
a neighborhood $N_{x}\subseteq U_{i}$ of $x$ such that, if $gN_{x}\cap F\neq\emptyset$, then $gN_{x}$ is a subset of some $U_{j}$. Since ${\overline W_{k}}$ is
compact, we can replace this $W_{k}$ by a finite number of neighborhoods like $N_{x}$  and the corresponding open sets $gN_{x}$ for $g\in M_{i}$. In this way, we
obtain a refinement of $\mathcal{U}$ which will be denoted again with $\mathcal{U}$, and shall be used in the sequel.

This refinement remains locally finite and, in addition, has the required property, because if $gN_{x}\cap F\neq\emptyset$, then $gU_{i}\cap F\neq\emptyset$, from
which follows that $g\in M_{i}$, hence $gN_{x}$ is a member of our refinement. It is easily seen that this property passes to the family $\{
St(x,\mathcal{U})\,|\,x\in F\}$, because if $gx\in
\overline{F}$, then $gSt(x,\mathcal{U})=St(gx,\mathcal{U})$.\\
\\
\textbf{3.4.} Next, using the covering $\mathcal{U}$ of ${\overline F}$ defined in 3.3, we consider the invariant covering $\mathcal{C}=\{gU_{i}\,|\,U_{i}\in
\mathcal{U},\,g\in G\}$ of $X$. We show that it is locally finite: For $x\in X$ there exists $h\in G$ such that $hx\in F$. Since $F$ is open, there exists an open and
relatively compact neighborhood $N\subseteq F$ of $hx$ that intersects finitely many members of $\mathcal{U}$. Then, the neighborhood $h^{-1}N$ of $x$ intersects
finitely many members of $\mathcal{C}$, because by 3.3, if $gU_{i}\cap N\neq\emptyset$, then
$gU_{i}\in \mathcal{U}$.\\
\\
\textbf{3.5.} To each covering like $\mathcal{C}$ corresponds a polyhedron $X_{\mathcal{C}}$, namely the nerve of the covering $\mathcal{C}$ with the CW-topology. A
subordinated partition of unity $\Phi _{\mathcal{C}}=\{\varphi _{U}\,|\,U\in \mathcal{C}\}$ determines a \textit{canonical map} $p_{\mathcal{C}}:X\rightarrow
X_{\mathcal{C}}$ with the property that $p_{\mathcal{C}}$ maps a point $x\in X$ to the point of $X_{\mathcal{C}}$ whose barycentric coordinate corresponding to the
vertex $U$ equals to $\varphi _{U}(x)$.

Since $\mathcal{C}$ is invariant, the initial properly discontinuous action $(G,X)$ induces a natural action $(G,X_{\mathcal{C}})$ defined as follows: For $g\in G$
and $U$ a vertex of $\mathcal{C}$ we let $(g,U)\mapsto gU$ and we extend the action map by linearity. This
action is properly discontinuous as it is easily verified.\\
\\
\textbf{3.6.} The construction of the desired inverse system of properly discontinuous actions on polyhedra will be based in the proof of the following theorem (cf.
\cite{Alder}, see also \cite[Th. 7 and Cor. 5, pp. 84-85]{Mar1}).\\
\\
\textbf{Theorem} \textit{Every connected, locally compact and paracompact space is the inverse limit of polyhedra.}\\
\\
For the convenience of the reader, we outline the proof: Let $X$ be a connected, locally compact paracompact space and {\Large$\mathcal{F}$} be the family of all
coverings of $X$  admitting a subordinated partition of unity. For every $\mathcal{D}\in$ {\Large$\mathcal{F}$} we choose a locally finite partition of unity $\Phi
_{\mathcal{D}}$ subordinated to $\mathcal{D}$. Let $X_{\mathcal{D}}$ be the nerve of $\mathcal{D}$ with the CW-topology. Let $ \Lambda $ be the set of all finite
subsets $\lambda =(\mathcal{D}_{1},\ldots ,\mathcal{D}_{n})$ of {\Large$\mathcal{F}$} ordered by inclusion. We denote by $X_{\lambda }$ the nerve of the covering
$$
\mathcal{D}_{1}\wedge \ldots \wedge \mathcal{D}_{n}=\{ V_{1}\cap \ldots \cap V_{n}\,|\,(V_{1},\ldots ,V_{n})\in \mathcal{D}_{1}\times \ldots \times \mathcal{D}_{n}\}.
$$
If $ \lambda \leq \mu =\{\mathcal{D}_{1},\ldots ,\mathcal{D}_{n},\ldots ,\mathcal{D}_{l}\}$, let $p_{\lambda \mu }:X_{\mu }\rightarrow X_{\lambda }$ be the simplicial
map which maps the vertex $(V_{1},\ldots ,V_{n},\ldots ,V_{l})$ of the nerve $\mathcal{D}_{1}\wedge \ldots \wedge \mathcal{D}_{n}\wedge \ldots \wedge \mathcal{D}_{l}$
to the vertex $ (V_{1},\ldots ,V_{n})$ of the nerve of $\mathcal{D}_{1}\wedge \ldots \wedge \mathcal{D}_{n}$.

As it is shown in \cite{Alder}, the family
$$
\Phi_{\mathcal{D}_{1}\wedge \ldots \wedge \mathcal{D}_{n}}=\{\varphi _{(V_{1},\ldots ,V_{n})}\,|\,(V_{1},\ldots ,V_{n})\in \mathcal{D}_{1}\times \ldots \times
\mathcal{D}_{n}\},
$$
where $\varphi _{(V_{1},\ldots ,V_{n})}=\varphi _{V_{1}}\cdot \ldots \cdot \varphi _{V_{n}}$, is a partition of unity subordinated to the covering
$\mathcal{D}_{1}\wedge \ldots \wedge \mathcal{D}_{n}$. Using this, for $\lambda =(\mathcal{D}_{1},\ldots ,\mathcal{D}_{n})$ we define the canonical map $p_{\lambda
}:X\rightarrow X_{\lambda }$ as in 3.5.

In order to obtain a polyhedral resolution of $X$ (cf. 1.4), a slight modification of the above construction is needed:

We replace the previous inverse system $(X_{\lambda },p_{\lambda \mu },\Lambda )$ by a larger system $(Y_{r},q_{rs},S)$ defined as follows: For $ \lambda \in \Lambda
$ let $\mathcal{V}_{\lambda }$ be a neighborhood basis of the closure of $p_{\lambda }(X)$ in $X_{\lambda }$, and
\begin{center}
$S=\{r=(\lambda ,V)\,|\,\lambda \in \Lambda\,\,\mbox{and}\,\, V\in \mathcal{V}_{\lambda}\}$.
\end{center}
Let $r\leq s=(\mu ,W)$ if $\lambda \leq \mu $ and $p_{\lambda \mu }(W)\subseteq V$. Moreover, letting $Y_{r}=V$, for $r\leq s$ we define the map
$q_{rs}:Y_{s}\rightarrow Y_{r}$  as the restriction of $p_{\lambda \mu}$ on $W$.

Taking into account the fact that {\Large$\mathcal{F}$} consists of \textit{all} coverings of $X$ admitting subordinated partitions of unity, it is proved
that $X=\varprojlim Y_{r}=\varprojlim X_{\lambda }$.\\
\\
\textbf{3.7.} If we replace {\Large$\mathcal{F}$} by {\Large$\mathcal{P}$}, the family of the coverings of $X$ of the form $\mathcal{C}=\{gU_{i}\,|\,U_{i}\in
\mathcal{U},\,g\in G\}$ defined in 3.4, and we repeat the previous steps, we obtain an inverse system, denoted (for simplicity) again by $(X_{\lambda },p_{\lambda\mu
},\Lambda )$. Since we use star coverings, we note that $St(x,\mathcal{D}_{1})\cap St(x,\mathcal{D}_{2})= St(x,\mathcal{D}_{1}\wedge \mathcal{D}_{2})$.\\
\\
\textbf{3.8.} If we restrict ourselves to the fundamental set ${\overline F}\subseteq X$, the coverings from {\Large$\mathcal{P}$} induce a family  of coverings on
${\overline F}$ defined by intersections of each one covering with ${\overline F}$. This family is cofinal to the corresponding one defined analogously via
{\Large$\mathcal{F}$} on ${\overline F}$. Since, {\Large$\mathcal{P}$} is not cofinal to {\Large$\mathcal{F}$} regarded as families of coverings of $X$, we shall
focus in the induced coverings of the fundamental set ${\overline F}$, where we may assume that the members of both families are the same. Note that, by \cite[ I,
Cor., p. 49]{Bour1}, ${\overline F}=\varprojlim p_{\lambda }({\overline F})$ holds, with respect to both {\Large$\mathcal{F}$} and {\Large$\mathcal{P}$}. Moreover,
with respect to {\Large$\mathcal{F}$}, we have $\varprojlim p_{\lambda }({\overline F})={\overline F}\subseteq X$, by the theorem in 3.6, while, with respect to
{\Large$\mathcal{P}$} and the notation from 3.7, ${\overline F}\subseteq\varprojlim X_{\lambda }$.

\section{The initial action as inverse limit of actions on polyhedra}
\begin{lemma}
Let $\mathcal{C}_{i}\in$ {\Large$\mathcal{P}$}. For the covering $\mathcal{C}_{1}\wedge \ldots \wedge \mathcal{C}_{n}$, there exists a subordinated partition of unity
$ \Phi_{\mathcal{C}_{1}\wedge \ldots \wedge \mathcal{C}_{n}}=\{\varphi _{(V_{1},\ldots ,V_{n})}\,|\,(V_{1},\ldots ,V_{n})\in \mathcal{C}_{1}\times \ldots \times
\mathcal{C}_{n}\}$ such that $\varphi_{(V_{1},\ldots ,V_{n})}=\varphi _{(gV_{1},\ldots ,gV_{n})}\circ g$, for every $g\in G$.
\end{lemma}
\begin{proof}
If the assertion is true for every single covering $\mathcal{C}$, then
$$
\varphi _{(gV_{1},\ldots ,gV_{n})}\circ g=[(\varphi_{V_{1}}\circ g^{-1})\cdot\ldots \cdot(\varphi_{V_{n}}\circ g^{-1})]\circ g=\varphi_{V_{1}} \cdot\ldots
\cdot\varphi_{V_{n}}=\varphi_{(V_{1},\ldots ,V_{n})}.
$$
\noindent So, it suffices to prove the assertion for a covering $\mathcal{C}=\{gU_{i}\,|\,U_{i}\in \mathcal{U},\,g\in G\}$ as in 3.4. We follow the usual construction
(cf. \cite[Th. 4.2, p. 170]{Dug1}): We choose locally finite coverings $\{V_{i}\,|\,i\in I\}$ and $\{W_{i}\,|\,i\in I\}$ of the open fundamental set $F$ such that
$\overline{W_{i}}\subset V_{i}\subset \overline{V_{i}}\subset U_{i}$ for every $i\in I$. We can apply Uryshon's Theorem in order to find continuous maps $f_{U_{i}}:X
\rightarrow [0,1]$ which are identically 1 on $\overline{W_{i}}$ and vanish on $X\setminus V_{i}$. We set $f_{gU_{i}}=f_{U_{i}}\circ g^{-1}$ for every $g\in G$. Since
the covering $\{g\overline{W_{i}}\,|\,i\in I,\,g\in G\}$ is locally finite, it follows that for each $x\in X$ at least one and at most finitely many $f_{gU_{i}}$ are
not zero, therefore $\sum f_{gU_{i}}$ is a well-defined continuous real-valued map on $X$ and is never zero. So, we can define the required partition of unity by
setting
$$
\varphi_{gU_{i}}(y)=\frac{f_{gU_{i}}(y)}{\sum f_{gU_{i}}(y)}.
$$
Since $x\in U_{i}$ iff $gx\in gU_{i}$, we have
$$
\varphi_{gU_{i}}(gx)=\frac{f_{gU_{i}}(gx)}{\sum f_{gU_{i}}(gx)}=\frac{f_{U_{i}}\circ g^{-1}(gx)}{\sum f_{U_{i}}\circ g^{-1}(gx)}=\frac{f_{U_{i}}(x)}{\sum
f_{U_{i}}(x)}=\varphi_{U_{i}}(x).
$$
\end{proof}

\begin{theorem}
With the notation from 3.7, $X$ is equivariantly homeomorphic to $\varprojlim X_{\lambda }$.
\end{theorem}
\begin{proof}
We recall that the actions $(G,X_{\lambda })$ defined in 3.5 induce an action on $\varprojlim X_{\lambda }$ as follows: Let $g\in G$ and $\mathbf{x}\in\varprojlim
X_{\lambda }$ with coordinates $p_{\lambda }(x)$. The coordinates of $g\mathbf{x}$ are $gp_{\lambda }(x)$, i.e., $p_{\lambda }(gx)=gp_{\lambda }(x)$. This action is
well defined since the maps $p_{\lambda \mu }:X_{\mu }\rightarrow X_{\lambda }$ have the property $p_{\lambda \mu }(gx_{\mu })=gp_{\lambda \mu }(x_{\mu })$ for every
$g\in G$ and $x_{\mu}\in X_{\mu}$, by the definition of the action on each $X_{\lambda }$.

An equivariant homeomorphism $f:X\rightarrow \varprojlim X_{\lambda }$ may defined in the following way: For $x\in X$ there exists some $g\in G$ such that $gx\in
{\overline F}$ (cf. 1.3). We let $f(x)$ to be the point with coordinates $p_{\lambda }(f(x))=g^{-1}(p_{\lambda }(gx))$. We will prove that $f$ is well defined: It
suffices to prove that, if $x\in X$ and $g\in G$ with $gx\in {\overline F}$, then $g^{-1}(p_{\lambda }(gx))$ is independent of the choice of $g$. Indeed, with the
notation from 3.5 and 3.6 let $x\in V_{1}\cap\ldots\cap V_{n}$. Then, by the definition of the actions on the polyhedra and the previous lemma, we have:
$$
\varphi_{(g^{-1}gV_{1},\ldots ,g^{-1}gV_{n})}(g^{-1}(p_{\lambda }(gx)))=\varphi_{(gV_{1},\ldots ,gV_{n})}(p_{\lambda }(gx))=\varphi_{(V_{1},\ldots ,V_{n})}(p_{\lambda
}(x)).
$$
Using this and the fact that $f$ is the identity map on the open fundamental set $F$, we can first verify that $f$ is equivariant and then a homeomorphism.
\end{proof}

\section{The embedding of the action in $\protect\varprojlim\varepsilon X_{\lambda}$ and its basic properties}
\begin{theorem}
The space $\varprojlim \varepsilon X_{\lambda }$ is a zero dimensional compactification of $X$. Moreover, $G$ acts on $\varprojlim \varepsilon X_{\lambda }$ and
$(G,X)$ is equivariantly embedded in $(G,\varprojlim \varepsilon X_{\lambda })$.
\end{theorem}
\begin{proof}
The simplicial maps $p_{\lambda \mu }:X_{\mu }\rightarrow X_{\lambda }$ are proper surjections. Hence, by Proposition 1.1.2, they have unique extensions $\varepsilon
p_{\lambda \mu }:\varepsilon X_{\mu }\rightarrow \varepsilon X_{\lambda }$ that map the space of the ends of $X_{\mu }$ onto that of $X_{\lambda }$. Furthermore,
$\varepsilon p_{\lambda \lambda }:\varepsilon X_{\lambda }\rightarrow \varepsilon X_{\lambda }$ is the identity map of $\varepsilon X_{\lambda }$, and for $\kappa
\leq \lambda $ and $\lambda \leq \mu$ we have $\varepsilon p_{\kappa \lambda }\circ \varepsilon p_{\lambda \mu }=\varepsilon p_{\kappa \mu }$. Hence, they define an
inverse limit, $\varprojlim \varepsilon X_{\lambda }$.

Using the fact that each $\varepsilon X_{\lambda }$ is a zero-dimensional compactification of $X_{\lambda }$ and applying Proposition 1.4.1, we see that $\varprojlim
\varepsilon X_{\lambda }$ is a zero-dimensional compactification of $\varprojlim X_{\lambda }$. By Corollary 1.1.4, the action of $G$ on $X_{\lambda}$ is extended to
an action on $\varepsilon X_{\lambda }$ such that the following equivariant diagram commutes
$$
\begin{CD}
(G,\varprojlim X_{\lambda })                  @>>\textrm{id}_{\textrm{G}}\times \textrm{h}> (G,\varprojlim \varepsilon X_{\lambda })\\
@V\textrm{id}_{\textrm{G}}\times \textrm{p}_{\lambda }VV        @VV\textrm{id}_{\textrm{G}}\times \varepsilon \textrm{p}_{\lambda }V\\
(G,X_{\lambda })                              @>>\textrm {id}_{\textrm{G}}\times \textrm{i}_{\lambda }> (G,\varepsilon X_{\lambda })
\end{CD}
$$
where $i_{\lambda }:X_{\lambda }\rightarrow \varepsilon X_{\lambda }$ are the inclusion maps, $p_{\lambda }:\varprojlim X_{\lambda }\rightarrow X_{\lambda }$ and
$\varepsilon p_{\lambda }:\varprojlim \varepsilon X_{\lambda }\rightarrow \varepsilon X_{\lambda }$ are projections, $id_{G}$ is the identity map of $G$ and
$h:\varprojlim X_{\lambda }\rightarrow \varprojlim \varepsilon X_{\lambda }$ is defined by setting $h_{\lambda}=i_{\lambda}$.

That $(G,X)$ embeds equivariantly in $(G,\varprojlim \varepsilon X_{\lambda })$ is an immediate consequence of Theorem 4.2 and the above diagram.
\end{proof}

\noindent \textbf{Remark} The example in Section 2 shows that $\varprojlim \varepsilon X_{\lambda }$ does not necessarily coincide with the end point compactification
$\varepsilon X$ of $X$. However,

\begin{proposition} If $X$ has finitely many ends, then $\varepsilon X=\varprojlim \varepsilon X_{\lambda }$.
\end{proposition}
\begin{proof}
Let $e_{i}$ for $i=1,2,\ldots ,n$ be the ends of $X$ and $V_{1},V_{2},\ldots ,V_{n}$ be open neighborhoods  of them in $\varepsilon X$, respectively, with disjoint
closures and boundaries lying in $X$. Then, the set $K=\varepsilon X\setminus\bigcup_{i=1}^{n}V_{i} $ is a compact subset of $X$. Let $e_{1}$, $e_{2}$ be two distinct
ends in $\varepsilon X$ with the same image in $\varprojlim \varepsilon X_{\lambda }$ via the projection map $p:\varepsilon X\rightarrow \varprojlim \varepsilon
X_{\lambda }$. Such a projection exists since, by Theorem 5.1,  $\varprojlim \varepsilon X_{\lambda }$ is a zero-dimensional compactification of $X$ and $\varepsilon
X$ is the maximal one. Therefore, $e_{1}$ and $e_{2}$ should have the same image under the composition map $\varepsilon p_{\lambda}\circ p$. With the notation from
3.4, this means that there is a subfamily of $\mathcal{C}$ with infinitely many members $gU_{i}$  with the property $g_{i}U_{i}\cap K\neq\emptyset$. Then, we can find
a sequence $\{x_{k}\}$ with $x_{k}\in g_{k}U_{k}\cap K$. Since $K$ is compact, we may assume that $\lim x_{k}=x\in K$, from which follows that the covering
$\mathcal{C}$ fails to be locally finite at $x$; a contradiction.
\end{proof}

The following proposition shows that, especially for equicontinuous actions, the sets $J(x)$, defined in 1.2, can be replaced by the \textit{limit sets}
$$
L(x)=\{ y\in X\,|\, \mbox{there exists a net}\, \{g_{i}\}\, \mbox{in}\,\, G\,\, \mbox{with}\,\, g_{i}\rightarrow\infty\,\, \mbox{and}\,\, \lim g_{i}x=y\},
$$
which are simpler to handle. The points of the sets $L(x)$ are the \textit{limit points} of the action.

\begin{proposition}
Let $(G,X)$ be an equicontinuous action of
a locally compact group $G$ on a locally compact space $X$. Then $J(x)=L(x)$ holds for every $x\in X$. Moreover, if the nets $\{x_{i}\}$ in $X$ and $\{g_{i}\}$ in $G$
are such that $\lim x_{i}=x$, $g_{i}\rightarrow \infty$ and $\lim g_{i}x_{i}=y\in J(x)$, then $\lim g_{i}x=y\in L(x)$.
\end{proposition}
\begin{proof}
Since $(G,X)$ is equicontinuous, for every entourage $U$ there exists an entourage $V$ such that for $y\in X$
$$
(x,y)\in V\, \mbox{implies}\, (gx,gy)\in U\, \mbox{for every}\, g\in G.
$$
But $\lim x_{i}=x$, so we may assume that $(x,x_{i})\in V$, therefore $(g_{i}x,g_{i}x_{i})\in U$ and $(g_{i}x_{i},y)\in U$. So $(g_{i}x,y)\in U\circ U$, hence $\lim
g_{i}x=y$.
\end{proof}

A kind of converse of the previous proposition is the following

\begin{proposition}
Let $Y$ be a zero-dimensional
compactification of the locally compact and connected space $Z$. Let $(G,Y)$ be an action such that $Z$ is an invariant subspace of $Y$ and the action $(G,Z)$ is
proper. The restricted action $(G,Z)$ is equicontinuous with respect to the uniformity induced on $Z$ by that of $Y$ iff the following condition is satisfied: If
$z\in Z$ is such that there exist a net $\{z_{i}\}$ in $Z$ with $\lim z_{i}=z$, and a net $\{g_{i}\}$ in $G$ with $g_{i}\rightarrow\infty$ and $\lim g_{i}z_{i}=e\in
Y\setminus Z$, then $\lim g_{i}z=e$.
\end{proposition}
\begin{proof}
The necessity may proved by arguments analogous to those applied in the proof of the previous proposition.

For the sufficiency, note that if the action $(G,Z)$ is not equicontinuous at the point $z$, then there exists an entourage $U$ such that for every entourage $V$
there exist a point $z_{\Vs}\in Z$ and some $g_{\Vs}\in G$ such that
$$
(z,z_{\Vs})\in V\, \mbox{and}\,  (g_{\Vs}z,g_{\Vs}z_{\Vs})\notin U.
$$
Since the entourages of the uniformity may directed by setting $V_{1}\leq V_{2}$ if $V_{2}\subseteq V_{1}$, we may assume that $g_{\Vs}\rightarrow\infty$ and $\lim
z_{\Vs}=z$. By the compactness of $Y$, we may assume that the nets $\{g_{\Vs}z\}$ and $\{g_{\Vs}z_{\Vs}\}$ converge to different points of $Y\setminus Z$, a
contradiction to our hypothesis.
\end{proof}

For the formulation of the next basic theorem, we recall that an action $(G,X)$ embedded in an action $(G,Y)$, where $Y$ is a zero-dimensional compactification of
$X$, is \textit{indivisible} if whenever $\lim g_{i}y_{0}=e\in Y\setminus X$ for some $y_{0}\in X$ then $\lim g_{i}y=e$ for every $y\in X$.

\begin{theorem} Let $(G,\varprojlim \varepsilon X_{\lambda })$
be the action defined in Theorem 5.1 and $L$ be the set of the limit points of the action $(G,X)=(G,\varprojlim X_{\lambda })$ in $\varprojlim \varepsilon X_{\lambda
}$. Then, the action $(G,\varprojlim \varepsilon X_{\lambda }\setminus L)$ is

(a) \hspace{0.01cm} proper,

(b) \hspace{0.01cm} equicontinuous with respect to the uniformity induced on $\varprojlim \varepsilon X_{\lambda }\setminus L$ by that of $\varprojlim \varepsilon
X_{\lambda }$, and

(c) \hspace{0.01cm} indivisible as embedded in the action $(G,\varprojlim \varepsilon X_{\lambda })$.
\end{theorem}

For the proof we need the following

\begin{lemma}
Let $L_{\lambda }$ be the set of the limit points of the action $(G,X_{\lambda})$ in $\varepsilon X_{\lambda }$. Then $L=\bigcap_{\lambda }\varepsilon p_{\lambda
}^{-1}(L_{\lambda })$.
\end{lemma}
\begin{proof}
Let $w\in L$, $\lim g_{i}x=w$ for $x\in X=\varprojlim X_{\lambda }$, and $g_{i}\rightarrow\infty$. Then $\lim g_{i}\varepsilon p_{\lambda }(x)=\lim \varepsilon
p_{\lambda }(g_{i}x)=\varepsilon p_{\lambda }(w)$, from which follows that $\varepsilon p_{\lambda }(w)\in L_{\lambda }$, therefore $L\subseteq\bigcap_{\lambda
}\varepsilon p_{\lambda }^{-1}(L_{\lambda })$.

For the inverse inclusion, let $v\in \bigcap_{\lambda }\varepsilon p_{\lambda }^{-1}(L_{\lambda })$, that is $\varepsilon p_{\lambda }(v)\in L_{\lambda }$ for every
$\lambda\in\Lambda$. This means that for each $\lambda\in\Lambda$ there exist a net $\{ g_{i}^{\lambda }\}$ in $G$ with $g_{i}^{\lambda }\rightarrow\infty$ and
$x_{\lambda }\in  X_{\lambda }$ with $\lim g_{i}^{\lambda }x_{\lambda }=\varepsilon p_{\lambda }(v)$. Since the polyhedron $X_{\lambda }$ is a connected, locally
compact and locally connected space, it has the ``property Z'', therefore, by \cite[3.4]{Ab2}, the action $(G,X_{\lambda })$ is indivisible as a restriction of
$(G,\varepsilon X_{\lambda })$. So, we may assume that $x_{\lambda }=\varepsilon p_{\lambda }(x)$ for a fixed $x\in X$ and every $\lambda\in\Lambda$. By the
compactness of $\varprojlim \varepsilon X_{\lambda }$, we may assume that $\lim g_{i}^{\lambda }x=v^{\lambda }\in L$. So, we have
$$
\varepsilon p_{\lambda }(v)=\lim g_{i}^{\lambda }x_{\lambda }=\lim g_{i}^{\lambda }\varepsilon p_{\lambda }(x)=\lim \varepsilon p_{\lambda }(g_{i}^{\lambda
}x)=\varepsilon p_{\lambda }(v^{\lambda }).
$$
Let $\lim v^{\lambda }=u\in\varprojlim \varepsilon X_{\lambda }$. This $u$ is contained in $\varprojlim \varepsilon X_{\lambda }\setminus X$, because $v^{\lambda
}\in\varprojlim \varepsilon X_{\lambda }\setminus X$ which, by Theorem 5.1, is a compact set. But, for each $\kappa\in\Lambda$ and every $\lambda\in\Lambda$ with
$\kappa\leq\lambda$, by 1.4, we have
$$
\varepsilon p_{\kappa }(u)=\lim \varepsilon p_{\kappa }(v^{\lambda })=\lim \varepsilon p_{\kappa \lambda }\circ \varepsilon p_{\lambda }(v^{\lambda })=\lim
\varepsilon p_{\kappa \lambda }\circ \varepsilon p_{\lambda }(v)=\varepsilon p_{\kappa }(v),
$$
from which follows that $u=v$. Taking into account this, that $\lim g_{i}^{\lambda }x=v^{\lambda }$ and applying a diagonal procedure, we may find a net $\{ g_{j}\}$
in $G$ such that $\lim g_{j}x=v\in\varprojlim \varepsilon X_{\lambda }\setminus X$. The properness of the action $(G,X)$ implies that this net is divergent, and
therefore $v\in L$, as required.
\end{proof}

\begin{proof}[Proof of Theorem 5.5.] (a) Assume that $\{ g_{i}\}$ is a net in $G$ and  $x$, $x_{i}$ and $y$  are points in $\varprojlim \varepsilon X_{\lambda
}\setminus L$ such that $\lim x_{i}=x$ and $\lim g_{i}x_{i}=y$. By the previous lemma, $\varprojlim \varepsilon X_{\lambda }\setminus L=\bigcup_{\lambda } \varepsilon
p_{\lambda }^{-1}(\varepsilon X_{\lambda} \setminus L_{\lambda })$. So, there exist $\kappa$ and $\lambda$ such that $x\in \varepsilon p_{\kappa}^{-1}(\varepsilon
X_{\kappa} \setminus L_{\kappa})$ and $y\in \varepsilon p_{\lambda }^{-1}(\varepsilon X_{\lambda} \setminus L_{\lambda })$.

For an index $\mu$ with $\kappa\leq\mu$ and $\lambda\leq\mu$, we may assume that
$$
\varepsilon p_{\kappa\mu}^{-1}(\varepsilon X_{\kappa }\setminus L_{\kappa })\cup \varepsilon p_{\lambda\mu}^{-1}(\varepsilon X_{\lambda}\setminus L_{\lambda
})\subseteq \varepsilon X_{\mu }\setminus L_{\mu }.
$$
Indeed, note that if, e.g., $z\in \varepsilon p_{\kappa\mu}^{-1}(\varepsilon X_{\kappa }\setminus L_{\kappa })$ and $z\in L_{\mu}$, there exist a net $\{ h_{j}\}$ in
$G$ with $h_{j}\rightarrow\infty$ and some $x_{\mu }\in X_{\mu }$ with $\lim h_{j}x_{\mu }=z$, hence $\lim h_{j}\varepsilon p_{\kappa\mu}(x_{\mu })=\varepsilon
p_{\kappa\mu}(z)\in L_{\kappa }$; a contradiction.

By this, we may assume that the points $x$, $x_{i}$ and $y$ are contained in the open and invariant set $\varepsilon p_{\mu}^{-1}(\varepsilon X_{\mu} \setminus
L_{\mu})$. Since, $X_{\mu }$ is connected, locally compact and locally connected, it has the ``property Z", therefore the action $(G,\varepsilon X_{\mu} \setminus
L_{\mu})$ is proper \cite[4.7]{Ab2}, hence $J(\varepsilon p_{\mu }(x))=\emptyset$. From this and the fact that $\lim g_{i}\varepsilon p_{\mu }(x_{i})=\varepsilon
p_{\mu }(y)$, it follows that the net $\{ g_{i}\}$ can not be divergent. Hence, by 1.2, the action $(G,\varprojlim \varepsilon X_{\lambda }\setminus L)$ is proper.

(b) We shall use Proposition 5.4 for $Z=\varprojlim \varepsilon X_{\lambda }\setminus L$ and the notation there. Let $\lim g_{i}z=e_{1}$. For every
$\lambda\in\Lambda$ we have
$$
\lim \varepsilon p_{\lambda }(z_{i})=\varepsilon p_{\lambda }(z),\,\,\lim g_{i}\varepsilon p_{\lambda }(z_{i})=\varepsilon p_{\lambda }(e),\,\,\mbox{and}\,\,\lim
g_{i}\varepsilon p_{\lambda }(z)=\varepsilon p_{\lambda }(e_{1}).
$$
By (a), the action $(G,Z)$ is proper, hence $e,e_{1}\in L$, therefore, by Lemma 5.6, $\varepsilon p_{\lambda }(e),\varepsilon p_{\lambda }(e_{1})\in L_{\lambda }$.
From this and the indivisibility of the action $(G,\varepsilon X_{\lambda }\setminus L_{\lambda})$ (cf. \cite[3.4]{Ab2}), it follows that $\varepsilon p_{\lambda
}(e)=\varepsilon p_{\lambda }(e_{1})$ for every $\lambda\in\Lambda$, i.e., $e=e_{1}$, and the assertion follows.

(c) The proof follows by repeating the arguments in the proof of (b).
\end{proof}

\noindent\textbf{Remark} We note that $X\subseteq\varprojlim \varepsilon X_{\lambda }\setminus L$, because $(G,X)$ is proper and $X=\varprojlim X_{\lambda }\subseteq
\varprojlim \varepsilon X_{\lambda }$.

\section{The maximality of $\protect\varprojlim\varepsilon X_{\lambda}=\mu X$ and the cardinality of $L=\mu L$}
In this paragraph we prove the main results of the paper.

\begin{lemma}
Let $(X,\mathbf{D})$ be a uniform space, and
$(G,X)$ be an equicontinuous action. Then, there exists a finer uniformity $\mathbf{D^{*}}$ compatible with the topology of $X$ such that $G$ acts on $X$ by
pseudoisometries with respect to the pseudometrics generating $\mathbf{D^{*}}$.
\end{lemma}
\begin{proof}
Let $\{d_{i},\, i\in I\}$ be a saturated family of bounded pseudome\-trics on $X$ which generates $\mathbf{D}$ \cite[II, Th. 1, p. 142]{Bour1}. We obtain a
pseudometric $d_{i}^{*}$ on $X$ such that every $h\in G$ acts on $X$ as a $d_{i}^{*}$-pseudoisometry, by letting $d_{i}^{*}(x,y)=\sup_{g\in G} d_{i}(gx,gy)$. Let
$\mathbf{D^{*}}$ denotes the uniformity generated by the family $\{ d_{i}^{*}\,|\,i \in I \}$. The topologies $\tau$, $\tau^{*}$ induced on $X$ by $\mathbf{D}$ and
$\mathbf{D^{*}}$, respectively, coincide: Since $d_{i}^{*}(x,y)\geq d_{i}(x,y)$, we have $\mathbf{D} \subseteq \mathbf{D^{*}}$  and $\tau \subseteq \tau^{*}$.
Conversely, if $U_{x}^{*} = \bigcap_{k=1}^{n} S_{k} (x,\epsilon)$ is a neighborhood of $x$ in $\tau^{*}$, where $S_{k}(x,\epsilon)$ denotes a $d^{*}_{i_{k}}$-ball of
radius $\epsilon$, centered at $x$, then the equicontinuity of $G$ implies the existence of a neighborhood $U_{x}$ of $x$ in $\tau$, such that $U_{x} \subseteq
U_{x}^{*}$.
\end{proof}

\begin{theorem}
The compactification $\varprojlim \varepsilon X_{\lambda }=\mu X$ is maximal among the zero-dimensional compactifications of $X$ satisfying simultaneously the
following properties:

(a) \hspace{0.01cm} The initial action $(G,X)$ is extended to an action $(G,\mu X)$.

(b) \hspace{0.01cm} The action $(G,\mu X\setminus\mu L)$, where $\mu L$ is the set of the limit points of the orbits of the initial action $(G,X)$ in $\mu X$, is
proper, equicontinuous with respect to the uniformity induced on $\mu X\setminus\mu L$ by that of $\mu X$, and indivisible.
\end{theorem}
\begin{proof}
By Theorem 5.5, the zero-dimensional compactification $\mu X$ of $X$ satisfies the properties (a) and (b). So, it remains to prove the maximality of $\mu X$: Suppose
that $Y$ is a zero-dimensional compactification of $X$ also satisfying these properties, such that $q:Y\rightarrow\varprojlim\varepsilon X_{\lambda }$ is a surjection
extending the identity map of $X$ (cf. 1.1). We have to show that $q$ is bijective.

\vspace{0.3cm}

\textit{Claim 1 : The restriction of $q$ on the set $L_{Y}$ of the limit points of the action $(G,X)$ in $Y$ is a bijection.}

Let $L_{Y}$ be the set of the limit points of the action $(G,X)$ in $Y$, and  $c_{1},\,c_{2}\in L_{Y}$ be two distinct points such that $q(c_{1})=q(c_{2})$. Then,
there are open neighborhoods $V_{1}$ and $V_{2}$ of $c_{1}$ and $c_{2}$, respectively, in $Y$ with disjoint closures.

Due to the indivisibility of the action $(G,Y\setminus L_{Y})$, we may assume that $\lim g_{j}x=c_{1}$ and $\lim h_{j}x=c_{2}$ with $g_{j}x\in V_{1}$ and $h_{j}x\in
V_{2}$. If there exists a covering $\mathcal{C}=\{gU_{i}\,|\,U_{i}\in \mathcal{U},\,g\in G\}$, as in 3.7, such that the members of it containing $g_{j}x$,
respectively $h_{j}x$, are pairwise disjoint, then, it is easily seen, that $\lim \varepsilon p_{\lambda }(q(g_{j}x))=\lim \varepsilon p_{\lambda }(g_{j}x)\neq \lim
\varepsilon p_{\lambda }(q(h_{j}x))$, a contradiction to the assumption that $q(c_{1})=q(c_{2})$. Therefore, there exist cofinal families $\{A_{j}\}$ and $\{B_{j}\}$
of members of $\mathcal{C}$ such that $g_{j}x\in A_{j}$, $h_{j}x\in B_{j}$ and $A_{j}\cap B_{j}\neq\emptyset$.

Since we consider star coverings, using refinements if it is necessary, we may assume that $A_{j}\cup B_{j}=f_{j}U_{j}$ is a member of our covering intersecting both
$V_{1}$ and $V_{2}$. Then $g_{j}x\in f_{j}U_{j}$, hence $x\in g_{j}^{-1}f_{j}U_{j}\in \mathcal{C}$. Since $\mathcal{C}$ is a locally finite covering, passing if
necessary to a subnet, we may assume that $x\in g_{j}^{-1}f_{j}U_{j}=gU_{r}$ for suitable $g$ and $r$. It follows that $h_{j}x=g_{j}gx_{j}$, where $x_{j}\in U_{r}$.
Since $U_{r}$ is relatively compact in $X$, we may assume that $\lim x_{j}=y\in X$. Thus, $c_{2}\in J(y)$ with respect to the action $(G,Y)$, because $\lim
h_{j}x=c_{2}$. From this and the assumption that the action $(G,Y\setminus L_{Y})$ is equicontinuous, taking into account Proposition 5.4, we conclude that $\lim
g_{j}gy=c_{2}$. This contradicts the fact that $gy\in X$, $\lim g_{j}x=c_{1}$ and the action $(G,Y\setminus L_{Y})$ is indivisible.

\vspace{0.3cm}

\textit{Claim 2 : The restriction of $q$ on the set $Y\setminus L_{Y}$ is also a bijection.}

Since $q$ is, by definition, the identity map on $X$, we have to show that it is bijective on $Y\setminus (L_{Y}\cup X)$. The action $(G,Y\setminus L_{Y})$ is
equicontinuous hence, by Lemma 6.1, we may assume that $G$ acts by pseudoisometries. So, we are allowed to assume that the invariant covering $\mathcal{C}$ consists
of open sets leading to invariant entourages.

Let $b_{1},b_{2}\in Y\setminus (L_{Y}\cup X)$ be two distinct points such that $q(b_{1})=q(b_{2})$, and $V=\{ (x,y)\in (Y\setminus L_{Y})\times (Y\setminus
L_{Y})\,|\,d_{k}(x,y)<\epsilon ,\,k=1,2,\ldots ,n\}$ be an entourage such that $(b_{1},b_{2})\notin V$. Moreover, we may assume that $\mathcal{C}$ consists of open
sets leading to invariant entourages of the form
$$
W=\{ (x,y)\in (Y\setminus L_{Y})\times (Y\setminus L_{Y})\,|\,d_{k}(x,y)<\epsilon /2 ,\,k=1,2,\ldots ,n,n+1,\ldots,m\}.
$$
As in the proof of Claim 1 and the notation there, since $b_{1},b_{2}\notin X$, we can find families $\{A_{j}\}$ and $\{B_{j}\}$ of members of $\mathcal{C}$ and
$x_{j}\in A_{j}$, $y_{j}\in B_{j}$ such that $\lim x_{j}=b_{1}$, $\lim y_{j}=b_{2}$ and $A_{j}\cup B_{j}$ is a member of our covering. From this and the specific
choice of the entourages $V$ and $W$, it follows that $(b_{1},b_{2})\in V$; a contradiction.
\end{proof}

\noindent\textbf{Corollary} If $X$ has the ``property Z", then $\mu X=\varepsilon X$.
\begin{proof}
We have to show that if $(G,X)$ is a properly discontinuous action, then the properties (a) and (b) of the previous theorem are satisfied for $\varepsilon X$, the
maximal zero-dimensional compactification of $X$, instead of $\mu X$. This follows from  Corollary 1.1.4, the already mentioned results of \cite{Ab2} in the
introduction (cf. \cite[4.7 and 3.7]{Ab2}), and by Proposition 5.4 for $Y=\varepsilon X$ and $Z=\varepsilon X\setminus \varepsilon L$.
\end{proof}

\noindent\textbf{Example} In our counterexample the set $\mu L$ consists of the two end points of the diagonal and the zero-dimensional compactification $\mu X$ may
obtained as a quotient space of the half-open Alexandroff square by identifying on the one hand the points $\{ (x,1)\,|\,x\in (0,1]\}$, and on the other hand the
points $\{ (x,0)\,|\,x\in [0,1)\}$.

\begin{theorem} The space $\mu L$ of the limit points of the
action $(G,X)$ in $\mu X$ either consists of at most two points or it is a perfect compact set. In the case where the group $G$ is abelian, $\mu L$ has at most two
points.
\end{theorem}
\begin{proof}
Let $\mu L$ has infinitely many points. We have to show that for every point $e$ of it and every neighborhood $V=\varepsilon p_{\lambda }^{-1}(U)$ of $e$ (cf.
Proposition 1.4.1), we can find a point $e_{1}\in\mu L$ with $e_{1}\in\overline{V}$ and $e_{1}\neq e$: By \cite[4.2]{Ab2}, the action $(G,\varepsilon X_{\lambda
}\setminus L_{\lambda })$ is proper. In the case under consideration, $L_{\lambda}$ is a perfect compact set, by \cite[4.11, Satz D, 4]{Ab2}. Since, by Lemma 5.6,
$\varepsilon p_{\lambda }(e)\in L_{\lambda }$, we can find $e_{\lambda }\in L_{\lambda }$ with $e_{\lambda }\in U$ and $e_{\lambda }\neq \varepsilon p_{\lambda }(e)$.
By the indivisibility of the action $(G,X_{\lambda })$ (cf. \cite[3.4]{Ab2}), there is a net $\{ g_{i}\}$ in $G$ with $g_{i}\rightarrow\infty$ and $\lim
g_{i}\varepsilon p_{\lambda }(x)=e_{\lambda }$ for every $x\in X$. Since $\lim \varepsilon p_{\lambda }(g_{i}x)=e_{\lambda }$ and $\mu X$ is compact, fixing an $x\in
X$, we may assume that $g_{i}x\in V$ and $\lim g_{i}x=e_{1}\in\mu L$. Therefore $e_{1}\in\overline{V}$. Since $\varepsilon p_{\lambda }(e_{1})=\lim \varepsilon
p_{\lambda }(g_{i}x)=e_{\lambda }\neq \varepsilon p_{\lambda }(e)$, we have $e_{1}\neq e$.

Now, assume that $\mu L$ has finitely many points. Since, by Lemma 5.6, $L=\bigcap_{\lambda }\varepsilon p_{\lambda }^{-1}(L_{\lambda })$, the set $\mu L$ is the
inverse limit of the inverse system $(L_{\lambda },\varepsilon p_{\lambda\mu }, \Lambda )$. By Theorem 1.2.1, every $L_{\lambda}$ consists of at most two points. From
this and the fact that every simplicial map $p_{\lambda \mu }:X_{\mu }\rightarrow X_{\lambda }$ is defined by deleting the last coordinates (cf. 3.6), we conclude
that $\mu X$ has also at most two points.

If the acting group is abelian, then the action $(G,X_{\lambda})$ fulfills the assumptions of the Theorem 1.17 of \cite{Lam1}, therefore every $L_{\lambda}$ consists
of one or two points. From this and using the same arguments as before, we see that $\mu L$ consists of at most two points.
\end{proof}

\section{An application}
In this section we apply our main results to show that the already known necessary condition for the existence of a proper action of a non-compact group on a locally
compact and connected space with the ``property Z" (cf. Theorem 1.2.1) remains also necessary in a broad class of actions, containing the properly discontinuous ones,
on spaces that do not have the ``property Z".

\begin{theorem} Let $X$  be a locally compact, connected and
paracompact space, and $G$ be a non-compact group acting properly on $X$ such that either $G_{0}$, the connected component of  the neutral element of $G$, is
non-compact, or $G_{0}$ is compact and $G/G_{0}$ contains an infinite discrete subgroup. Then $X$ has

(a) \hspace{0.01cm} at most two or infinitely many ends, and

(b) \hspace{0.01cm} at most two ends, if $G_{0}$ is not compact.
\end{theorem}
\begin{proof}
We begin with the proof of (b) and we shall restrict ourselves in the proof of (a) only in the case where $G_{0}$ is compact.

(b) If $G_{0}$ is non-compact, we consider the restricted action $(G_{0},X)$. By Iwasawa's Decomposition Theorem, $G_{0}$ contains a closed subgroup isomorphic to
$\mathbb{R}$, therefore it contains a closed subgroup isomorphic to $\mathbb{Z}$, the additive group of the integers. The restricted action $(\mathbb{R} ,X)$ is
proper, therefore the action $(\mathbb{Z} ,X)$ is properly discontinuous.

Since the space of the ends of $X$ is totally disconnected, every end is a fixed point for the action $(G_{0},\varepsilon X)$, therefore for the restricted action
$(\mathbb{Z} ,\varepsilon X)$ too. Since the projection $p:\varepsilon X\rightarrow\mu X$ is equivariant, every point of $\mu X\setminus X$ is a fixed point for the
action $(\mathbb{Z} ,\mu X)$, where $\mu X$ is the zero-dimensional compactification of $X$ that corresponds to the action $(\mathbb{Z} ,X)$ by Theorem 6.2. From this
and Theorem 6.3, there exist at most two limit points for the action $(\mathbb{Z} ,\mu X\setminus \mu L)$. The set $\mu X\setminus X$ cannot have any other point
except these limit points, because by Theorem 5.5, the action $(\mathbb{Z} ,\mu X\setminus\mu L)$ is properly discontinuous, therefore has compact isotropy groups.

We claim that in this case $\varepsilon X=\mu X$ holds, which implies (b). To this end, we have to prove that $p$ is injective. In order to be able to repeat the
arguments in the proof of Theorem 6.2, Claim 2 replacing $Y$ by $\varepsilon X$,  we need the following

\vspace{0.3cm}

\textit{Claim : The action $(\mathbb{R},X)$ is equicontinuous with respect to the uniformity induced on $X$ by that of $\varepsilon X$.}

We shall use Proposition 5.4. Let $x\in X$ and $\lim x_{i}=x$ for $x_{i}\in X$. To arrive at a contradiction, assume that there exists a net $\{ t_{i}\}$ in
$\mathbb{R}$ with $t_{i}\rightarrow +\infty$ and $\lim t_{i}x=e_{1}\in \varepsilon X\setminus X$, while $\lim t_{i}x_{i}=e_{2}\in \varepsilon X\setminus X$, where
$e_{1}\neq e_{2}$. Let $U$ and $V_{1}$ be disjoint neighborhoods in $\varepsilon X$ of $x$ and $e_{1}$, respectively, with boundaries in $X$. Then, there exists
$t_{0}$ such that $tx\in V_{1}$ for every $t\geq t_{0}$, because otherwise, by the connectedness of the orbits, we can find a net $\{ r_{i}x\}$ in the boundary of
$V_{1}$ with $\lim r_{i}x=y\in X$ and $r_{i}\rightarrow +\infty$; this is not possible, because the action $(\mathbb{R} ,X)$ is proper, hence $L(x)\subseteq
J(x)=\emptyset$ for every $x\in X$ (cf. 1.2 and Section 5). So, we can find a neighborhood $V_{2}$ in $\varepsilon X$ of $e_{2}$ with boundary in $X$, disjoint from
$U$ such that $tx\notin\overline{V_{2}}$ for every $t\geq 0$. Since $\lim x_{i}=x\in U$ and $\lim t_{i}x_{i}=e_{2}\in V_{2}$ there exists a net $\{ s_{i}x_{i}\}$ in
the boundary of $V_{2}$ such that $\lim s_{i}x_{i}=z\in X$. As before, the net $\{s_{i}\}$ cannot be divergent, therefore we may assume that $\lim s_{i}=s\geq 0$.
Hence $z=sx\in\overline{V_{2}}$; a contradiction.

(a) We have to consider only the case where $G_{0}$ is compact and $G/G_{0}$ contains an infinite discrete subgroup. Since $X$ is connected and $\sigma$-compact, the
orbit space $X\backslash G_{0}$ of the action $(G_{0},X)$ is connected and $\sigma$-compact, therefore paracompact. The group $G/G_{0}$ acts on $X\backslash G_{0}$,
by letting
$$
(gG_{0},G_{0}(x))\mapsto G_{0}(gx),\,\,\mbox{for every}\,\,g\in G\,\,\mbox{and}\,\,x\in X.
$$
This action is well defined, since $G_{0}$ is a normal subgroup of $G$. Moreover, it is proper: Since the initial action is proper, $G$ is locally compact, therefore,
by 1.2, there exist compact neighborhoods $U_{x}$, $U_{y}$ in $X$ of $x$ and $y$, respectively, such that the set
$$
G(U_{x},U_{y})=\{ g\in G\,|\,(gU_{x})\cap U_{y}\neq\emptyset\}
$$
is relatively compact in $G$. Then $W_{1}=\{G_{0}(z)\,|\,z\in U_{x}\}$ and $W_{2}=\{G_{0}(z)\,|\,z\in U_{y}\}$ are compact neighborhoods of the points $G_{0}(x)$,
$G_{0}(y)$ in $X\backslash G_{0}$, respectively. The set
$$
(G/G_{0})(W_{1},W_{2})=\{ gG_{0}\in G/G_{0}\,|\,(gG_{0}W_{1})\cap W_{2}\neq\emptyset\}
$$
is relatively compact in $G/G_{0}$. Indeed, let $\{g_{i}G_{0}\}$ be a net in $(G/G_{0})(W_{1},W_{2})$. Then, there exist $h_{i}$, $q_{i}\in G_{0}$ and $x_{i}\in
U_{x}$, $y_{i}\in U_{y}$ such that $g_{i}h_{i}x_{i}=q_{i}y_{i}$, i.e., $q_{i}^{-1}g_{i}h_{i}\in G(U_{x},U_{y})$. Therefore $g_{i}\in G_{0}\cdot G(U_{x},U_{y})\cdot
G_{0}$ which is a relatively compact subset of $G$. This means that $\{g_{i}G_{0}\}\rightarrow\infty$ is not possible. Hence $(G/G_{0})(W_{1},W_{2})$ is relatively
compact. Therefore, every non-compact discrete subgroup $F$ of $G/G_{0}$ acts properly discontinuously on $X\backslash G_{0}$, which is a locally compact, connected
and paracompact space.

So, we can apply our results for the action $(F,X\backslash G_{0})$. The surjective map $q:X\rightarrow X\backslash G_{0}$, with $q(x)=G_{0}(x)$ is proper, because
$G_{0}$ is compact, therefore, by Proposition 1.1.2, it has a unique extension $\varepsilon q:\varepsilon X\rightarrow \varepsilon (X\backslash G_{0})$ that maps the
ends of $X$ onto that of $X\backslash G_{0}$. So, this map relates the ends of $X$ with those of $X\backslash G_{0}$.

\vspace{0.3cm}

\textit{Claim : The restriction of the map $\varepsilon q$ on the set of the ends of $X$ is a bijection.}

Since $G_{0}$ is connected, as before, the ends of $X$ are fixed points for the action $(G_{0},\varepsilon X)$. The map $\varepsilon q$ is equivariant, therefore the
ends of $X\backslash G_{0}$ are also fixed points with respect to the action $(G/G_{0},\varepsilon (X\backslash G_{0}))$. Since every end of $X$ is a $G_{0}$-orbit,
the assertion follows.

\vspace{0.3cm}

If $X$ has infinitely many ends there is nothing to prove. If $X$ has finitely many ends, let $\mu (X\backslash G_{0})$ be the zero-dimensional compactification of
$X\backslash G_{0}$ that corresponds to the action $(F,X\backslash G_{0})$ by Theorem 6.2. According to Proposition 5.2, we have that $\mu (X\backslash
G_{0})=\varepsilon (X\backslash G_{0})$, and by Theorem 6.3, the set $L^{*}$ of the limit points of the action $(F,X\backslash G_{0})$ consists of at most two points.
There are no other ends except those of $L^{*}$, because by Theorem 6.2(b), the non-compact group $F$ acts properly on $\varepsilon (X\backslash G_{0})\setminus
L^{*}$ which has finitely many points. This and the previous claim prove the theorem.
\end{proof}

\bibliographystyle{amsplain}

\end{document}